%% file: main.tex
\documentclass[oneside, english]{amsart}

\usepackage{geometry}
\usepackage{amsmath, amsthm}
\usepackage[hidelinks]{hyperref}
\usepackage{cleveref}
\usepackage{amssymb, stmaryrd, mathrsfs}
\usepackage{enumitem}

\theoremstyle{definition}
	\newtheorem{defi}{Definition}[section]
	\newtheorem{rema}[defi]{Remark}
	
\theoremstyle{plain}
	\newtheorem{theo}[defi]{Theorem}
	\newtheorem{prop}[defi]{Proposition}
	\newtheorem{lemm}[defi]{Lemma}
	\newtheorem{coro}[defi]{Corollary}
	
\newcommand{\sm}{\setminus}
\newcommand{\mr}[1]{\mathrm{#1}}
\newcommand{\mc}[1]{\mathcal{#1}}
\newcommand{\bb}[1]{\mathbb{#1}}
\newcommand{\val}[1]{\left| #1 \right|}

\DeclareMathOperator{\Th}{Th}
\DeclareMathOperator{\graph}{graph}
\DeclareMathOperator{\vox}{vox}

\renewcommand{\leq}{\leqslant}
\renewcommand{\subset}{\subseteq}
\setlist[enumerate]{label={\upshape(\roman*)}}

\title[Definably isometric cuboids over Hensel minimal fields]{%
	Definably isometric cuboids over\\
	Hensel minimal fields
}

\author{Marcelo Fortmann}

\address{%
	Lehrstuhl für Algebra und Zahlentheorie, 
	Mathematisches Institut, 
	Heinrich-Heine-Universit\"at \allowbreak D\"usseldorf, 
	Universit\"atsstraße 1, 40225 D\"usseldorf, 
	Germany
}

\email{marcelo.fortmann@hhu.de}

\def\myabstract{%
	We show that, over a valued field $K$ with $1$-h-minimal theory $\Th(K)$,
	two Cartesian products of balls must be equal up to translation and coordinate permutation
	in order for there to exist a definable surjective isometry between them. 
	For this purpose, we introduce a notion of relative dimension for cuboids in $K^n$
	that is computed combinatorially and relates to the notion of dimension for definable sets
	in Hensel minimal theories when $K$ has infinite residue field.
}

\def\mythanks{%
	The author would like to thank Immanuel Halupczok for giving them this question 
	and for many helpful discussions throughout all stages of writing this article.
}

\def\myai{%
	During the creation of this article,
	the author used a free plan of ChatGPT for broad mathematical discussion.
	The proofs of \Cref{lemm:dim-type_intersection_with_ball}
	and \Cref{lemm:isometry_descends_to_bijection} were suggested by ChatGPT.
	The entire article, including the aforementioned proofs, was written by the author.
}

\begin{document}

\begin{abstract}
\myabstract
\end{abstract}

\maketitle

\input{section_1}

\subsection{Acknowledgements}
\mythanks

\subsection{AI statement}
\myai

\input{section_2}
\input{section_3}

\bibliographystyle{plain}
\bibliography{lit}

\end{document}

%% file: section_1.tex
\section{Introduction}

Let $K$ be a valued field with infinite residue field $k$.
Then any bijection from $k$ to $k^2$ lifts to an isometric bijection
from $\mc O \times \mc M$ to $\mc O \times \mc O$,
where $\mc O$ and $\mc M$ denote the valuation ring and valuation ideal of $K$, respectively. 
Precisely when such a map is first-order definable is not clear.
In general, it is not clear when there exists a definable surjective isometry
between two Cartesian products of balls in $K$.
We show that, when $\Th(K)$ is $1$-h-minimal, 
two Cartesian products of balls in $K$ must be equal up to translation and coordinate permutation
in order to be definably isometric.
This question was asked by Immanuel Halupczok
and our presented answer will be used in an upcoming article 
by Cluckers, Halupczok and Vermeulen.

\subsection{Main results}

Let $K$ be a non-trivially valued field,
viewed as a structure in a language that contains 
the language of valued fields $\mc L_{\mr{val}} = \{+,\cdot,\mc O\}$.
Let $\left| \,\cdot\, \right|$ denote the valuation on $K$, multiplicatively written,
that is extended to $K^n$ by the sup norm.
We call a subset of $K^n$ definable if it is definable with any parameters from $K$. 
Our main result is the following.

\begin{theo}[\Cref{theo:infinite_residue_field_isometric_cuboids_equal_up_to}]
\label{theo:1.1}
Assume that $\Th(K)$ is $1$-h-minimal and let $P,P' \subset K^n$ be Cartesian products of balls in $K$.
If $P$ and $P'$ are definably isometric, then $P$ and $P'$ are equal 
up to translation and coordinate permutation.
\end{theo}

In particular, this holds true regardless of the residue characteristic of $K$.
Our notion of $1$-h-minimality therefore refers to the general case described in \cite[Theorem 2.2.8]{Cluckers_Halupczok_Rideau-Kikuchi_Vermeulen_2023}.
Cartesian products of balls in $K$ are the main object of study in this article 
and will be called cuboids throughout.
A precise definition of cuboids and isometries is given in Section 2.
For the purpose of proving \Cref{theo:1.1}, 
we introduce a notion of relative dimension $\dim(P/Q)$ for cuboids $P,Q \subset K^n$.
This notion of dimension is defined purely combinatorially 
is consistent with the notion of dimension for definable subsets of $K^n$
in the sense of \cite[Definition 5.3.1]{Cluckers_Halupczok_Rideau-Kikuchi_2022}.
More precisely, for a cuboid $P \subset K^n$ around $0$
and a ball $B \subset K^n$ around $0$,
we have the following characterization.

\begin{prop}[\Cref{prop:dim_for_cuboids_relates_to_dim_for_definable_sets}]
\label{prop:1.2}
Assume that $\Th(K)$ is $1$-h-minimal and that $K$ has infinite residue field.
Then $\dim(P/B)$ is the minimal dimension of a definable set $X \subset K^n$
such that $\vox_B(X)$ covers $\vox_B(P)$.
\end{prop}

Here, $\vox_B$ denotes the natural projection from $K^n$ to $K^n/B$.
This map is essential for our proof of \Cref{theo:1.1}
and will be used all throughout Section 3.
The proof of \Cref{theo:1.1} is divided into two parts. 
When $K$ has finite residue field, the claim of \Cref{theo:1.1} applies to all isometric cuboids, 
regardless of definability. The techniques used in this part are elementary.
When $K$ has infinite residue field, the key ingredient to the proof is the following result.

\begin{prop}[\Cref{prop:counting_balls}]
\label{prop:1.3}
Assume that $\Th(K)$ is $1$-h-minimal. 
For every definable set $X \subseteq K^n$ there exists some $N \in \bb N$
such that for every ball $B \subset K$ around $0$ and every finite subset $A \subset K$, one has
\begin{equation*}
\#(\vox_{B^n}(X) \cap \vox_{B^n}(A^n)) \leq N \cdot (\#\vox_B(A))^{\dim(X)}.
\end{equation*} 
\end{prop}

Here, $\vox_{B^n}$ and $\vox_{B}$ denote the natural projections from $K^n$ to $K^n/B^n$
and from $K$ to $K/B$, respectively. For a precise definition, see Section 3.

%% file: section_2.tex
\section{Relative dimension of cuboids}

Throughout, let $\pi_i \colon K^n \to K$ denote the projection to the $i$th coordinate
for $1 \leq i \leq n$.

\begin{defi}
A \emph{ball in $K^n$} is a proper and infinite subset $B$ of $K^n$
such that for all $x,x' \in B$ and $y \in K^n$, $\val{x-y} \leq \val{x-x'}$ implies $y \in B$.
\end{defi}

Note that this definition includes balls derived from cuts in the value group.
Further note that whether or not a given set is a ball depends on the ambient space.
For example, if $L/K$ is an extension of valued fields and $B$ is a ball in $K$,
then $B$ is not necessarily a ball in $L$.
Since the ambient space will always be apparent, we often omit to mention it.

\begin{defi}
A \emph{cuboid in $K^n$} is a Cartesian product of $n$ balls in $K$.
\end{defi}

Clearly, every ball in $K^n$ is a cuboid in $K^n$.
Like balls, any cuboid in $K^n$ that contains $0$
is a subgroup of $K^n$ and the non-empty intersection of two cuboids is another cuboid.

\begin{defi}
An \emph{isometry} between subsets $X,Y \subseteq K^n$ is a bijection $f \colon X \to Y$
such that 
\begin{equation*}
\val{f(x)-f(x')} = \val{x-x'}
\quad \text{for all $x,x' \in X$.}
\end{equation*}
We call $X$ and $Y$ \emph{isometric} if there exists an isometry between them.
Moreover,
we call $X$ and $Y$ \emph{definably isometric} if there exists a definable isometry between them.
\end{defi}

As shown in the introduction, whenever $K$ has infinite residue field,
there exists an isometry between $\mc O \times \mc M$ and $\mc O \times \mc O$,
where $\mc O$ and $\mc M$ denote the valuation ring and valuation ideal of $K$, respectively.
We will see that such an isometry cannot be definable
by showing that definable isometries between cuboids preserve a certain sense of size and shape.
To make this precise, we introduce the following notion of relative dimension for cuboids.

\begin{defi}
Let $P,Q \subset K^n$ be cuboids around $0$.
\begin{enumerate}
\item
The \emph{dimension type of $P$ over $Q$} is the tuple $t \in \{0,1\}^n$
with entries $t_i = 0$ if and only if $\pi_i(P) \subseteq \pi_i(Q)$
for $1 \leq i \leq n$. We denote it by $t(P/Q)$.
\item
The \emph{dimension of $P$ over $Q$} is the sum of entries of $t(P/Q)$.
We denote it by $\dim(P/Q)$.
\end{enumerate}
\end{defi}

Clearly, translation does not affect
whether two cuboids are definably isometric.
Hence, we may assume without loss of generality that 
all cuboids in discussion contain $0$.

We will see that definably isometric cuboids have the same dimension
over certain balls in the sense of the above definition.
In particular, it will follow from \Cref{prop:infinite_residue_field_isometric_cuboids_same_dim} 
that the cuboids $\mc O \times \mc M$ and $\mc O \times \mc O$ 
must have the same dimension over every ball in $K^2$ in order to be definably isometric. 
This fails over $\mc M \times \mc M$, 
as $\dim(\mc O \times \mc M/\mc M \times \mc M) = 1$ and $\dim(\mc O \times \mc O/\mc M \times \mc M) = 2$.

\begin{rema}
Let $P,Q \subset K^n$ be cuboids around $0$. Then $t(P/Q \cap P) = t(P/Q)$.
\end{rema}

\begin{lemm}\label{lemm:dim-type_intersection_with_ball}
Let $P \subset K^n$ be a cuboid around $0$ and let $B \subsetneq B' \subset K^n$ be strictly nested balls around $0$.
Then $t(P \cap B'/B) = t(P/B)$. 
\end{lemm}

\begin{proof}
Fix some index $i$.
Clearly, the inclusion $\pi_i(P) \subseteq \pi_i(B)$ implies $\pi_i(P \cap B') \subseteq \pi_i(B)$.
Now suppose $\pi_i(P) \nsubseteq \pi_i(B)$.
Since $P$ and $B'$ have non-empty intersection, the balls $\pi_i(P)$ and $\pi_i(B')$ are nested.
Hence, $\pi_i(P \cap B')$ equals $\pi_i(P)$ or $\pi_i(B')$
and in either case one has $\pi_i(P \cap B') \nsubseteq \pi_i(B)$.
\end{proof}

%% file: section_3.tex
\section{Definable isometries between cuboids}

In this section, we go about proving \Cref{theo:1.1}.
We treat $K$ having finite and infinite residue field separately.
For a ball $B \subset K^n$ around $0$, we often make use of the natural projection
\begin{equation*}
\vox_B \colon K^n \to K^n/B.
\end{equation*}
The name $\vox_B$ is short for \emph{voxelization} with voxels the size of $B$.

\begin{rema}\label{rema:dim_existence_of_definable_cover}
Let $P \subset K^n$ be a cuboid and $B \subset K^n$ a ball, each around $0$.
Then $t(P/B)$ is the unique tuple $t \in \{0,1\}^n$
with the least sum of entries such that
$\vox_B(\prod_{i=1}^n t_i K)$ covers $\vox_B(P)$.
Here, $0 \cdot K = \{0\}$ and $1 \cdot K = K$.
\end{rema}

\begin{lemm}\label{lemm:isometry_descends_to_bijection}
Let $P,P' \subset K^n$ be isometric cuboids and $B \subset K^n$ a ball, each around $0$.
Then 
\begin{enumerate}
\item\label{lemm:isometry_descends_to_bijection-1}
the cuboids $P \cap B$ and $P' \cap B$ are isometric, and
\item\label{lemm:isometry_descends_to_bijection-2}
there exists a bijection between $\vox_B(P)$ and $\vox_B(P')$.
\end{enumerate}
\end{lemm}

\begin{proof}
Let $f \colon P \to P'$ be an isometry. 
For \ref{lemm:isometry_descends_to_bijection-1}, see that $f(P \cap B) = P' \cap (f(0) + B)$.
Hence, the mapping $x \mapsto f(x) - f(0)$ gives an isometry from $P \cap B$ to $P' \cap B$.
For \ref{lemm:isometry_descends_to_bijection-2}, 
all $x,x' \in K^n$ satisfy $f(x)-f(x') \in B$ if and only if $x-x' \in B$.
Hence, $f$ descends to a well-defined bijection from $\vox_B(P)$ to $\vox_B(P')$.
\end{proof}

\subsection{Finite residue field}

When $K$ has finite residue field, we can use \Cref{lemm:isometry_descends_to_bijection}
to find a finite bijection between $\vox_B(P)$ and $\vox_B(P')$ 
for suitable isometric cuboids $P,P' \subset K^n$ and balls $B \subset K^n$.

\begin{rema}\label{rema:finite_residue_field_count_vox}
Assume that $K$ has finite residue field $k$. Let $P \subset K^n$ be a cuboid around $0$ and let
$B = \{ x \in K^n \mid \val x < \gamma \}$ and $B' = \{ x \in K^n \mid \val x \leq \gamma \}$
for some $\gamma$ in the value group of $K$.
If $P$ is contained in $B'$, then $\#\vox_B(P) = (\#k)^{\dim(P/B)}$.
\end{rema}

\begin{lemm}\label{lemm:finite_residue_field_isometric_cuboids_same_dimension}
Assume that $K$ has finite residue field $k$. Let $P,P' \subset K^n$ be isometric cuboids around $0$
and let $B = \{ x \in K^n \mid \val x < \gamma \}$ for some $\gamma$ in the value group of $K$.
Then $\dim(P/B) = \dim(P'/B)$.  
\end{lemm}

\begin{proof}
Let $B' = \{ x \in K^n \mid \val x \leq \gamma \}$.
By \Cref{lemm:isometry_descends_to_bijection},
there exists a bijection between $\vox_B(P \cap B')$ and $\vox_B(P' \cap B')$.
Hence, \Cref{rema:finite_residue_field_count_vox} yields $\dim(P \cap B'/B) = \dim(P' \cap B'/B)$.
Now apply \Cref{lemm:dim-type_intersection_with_ball}.
\end{proof}

\begin{prop}\label{prop:finite_residue_field_isometric_cuboids_equal_up_to}
Assume that $K$ has finite residue field and let $P,P' \subset K^n$ be cuboids.
If $P$ and $P'$ are isometric, 
then $P$ and $P'$ are equal up to translation and coordinate permutation.
\end{prop}

\begin{proof}
Assume that $P$ and $P'$ are isometric.
Since translation and coordinate permutation do not affect
whether $P$ and $P'$ are isometric, we may assume without loss of generality 
that $P$ and $P'$ each contain $0$ 
and that $\pi_{i}(P) \subseteq \pi_{i+1}(P)$ and $\pi_{i}(P') \subseteq \pi_{i+1}(P')$ for $1 \leq i < n$.
Suppose that $P$ and $P'$ are not equal.
Then there is a largest index $i$ such that $\pi_i(P)$ and $\pi_i(P')$ differ. 
Without loss, assume $\pi_i(P) \subset \pi_i(P')$. 
By choice of $i$, this implies $\pi_i(P) \subset \pi_{i+1}(P)$. 
Pick some $a \in \pi_i(P') \sm \pi_i(P)$ and let $I = \{ x \in K \mid \val x < \val a \}$. 
Then $\pi_i(P) \subseteq I \subset \pi_i(P')$, 
hence $\dim(P/I^n) = n-i$ and $\dim(P'/I^n) = n-(i-1)$. 
But this contradicts \Cref{lemm:finite_residue_field_isometric_cuboids_same_dimension}.
\end{proof}

\subsection{Infinite residue field}

When $K$ has infinite residue field, 
the previous strategy fails, since $\vox_B(P)$
is either a singleton or infinite for any cuboid $P \subset K^n$ and ball $B \subset K^n$.
To bypass this issue, we assume that the theory of $K$ is $1$-h-minimal.
This way, we are able to get finite approximations of $\vox_B(X)$
for any definable set $X \subseteq K^n$ and ball $B \subset K^n$ in the form of \Cref{prop:1.3}.
\begin{defi}
A map $f \colon X \to K^n$ on a subset $X \subseteq K^n$ is called \emph{$1$-Lipschitz} if
\begin{equation*}
\val{f(x)-f(x')} \leq \val{x-x'}
\quad \text{for all $x,x' \in X$.}
\end{equation*}
\end{defi}

For the rest of this section,
let $\pi_{\leq d} \colon K^n \to K^d$ denote the projection to the first $d$ coordinates
for $0 \leq d \leq n$.
The following theorem is a simple corollary of 
\cite[Theorem 3.3.3, Addendum 4]{Cluckers_Halupczok_Rideau-Kikuchi_Vermeulen_2023}.

\begin{theo}[$1$-Lipschitz cell decomposition]\label{theo:1-Lipschitz_cell_decomposition}
Assume that $\Th(K)$ is $1$-h-minimal.
Then for any definable set $X \subseteq K^n$,
there exists a finite partition of $X$ into definable subsets $C_1,\ldots,C_N \subseteq X$,
together with coordinate permutations $\sigma_1,\ldots,\sigma_N \colon K^n \to K^n$,
such that each $\sigma_i(C_i)$
is the graph of a definable $1$-Lipschitz map
\begin{equation*}
\pi_{\leq \dim(C_i)}(\sigma_i(C_i)) \to K^{n - \dim(C_i)}.
\end{equation*} 
\end{theo}

\begin{rema}\label{rema:1-Lipschitz_map_induces_bijection}
Let $X \subseteq K^d$ be a subset for some integer $d \leq n$.
Then any $1$-Lipschitz map $f \colon X \to K^{n-d}$ satisfies 
\begin{equation*}
\val{(x,f(x)) - (x',f(x'))} = \val{x-x'}
\quad \text{for all $x,x' \in X$.}
\end{equation*}
Hence, for any ball $B \subset K$ around $0$,
the projection $\pi_{\leq d} \colon \graph(f) \to X$ 
descends to a well-defined bijection $\bar \pi_{\leq d} \colon \vox_{B^n}(\graph(f)) \to \vox_{B^d}(X)$.
\end{rema}

With this, we are able to prove \Cref{prop:1.3}.

\begin{prop}[Counting balls]\label{prop:counting_balls}
Assume that $\Th(K)$ is $1$-h-minimal. 
For every definable set $X \subseteq K^n$ there exists some $N \in \bb N$
such that for every ball $B \subset K$ around $0$ and every finite subset $A \subset K$, one has
\begin{equation*}
\#(\vox_{B^n}(X) \cap \vox_{B^n}(A^n)) \leq N \cdot (\#\vox_B(A))^{\dim(X)}.
\end{equation*} 
\end{prop}

\begin{proof}
Let $C_1,\ldots,C_N$ and $\sigma_1,\ldots,\sigma_N$ 
be as in \Cref{theo:1-Lipschitz_cell_decomposition}
and fix some index $i \in \{1,\ldots,N\}$.
Let $B \subset K$ be a ball around $0$ and let $A \subset K$ be a finite subset.
Since $\sigma_i(A^n) = A^n$, one has
\begin{equation*}
\#(\vox_{B^n}(C_i) \cap \vox_{B^n}(A^n))
= \#(\vox_{B^n}(\sigma(C_i)) \cap \vox_{B^n}(A^n)).
\end{equation*}
Since $\sigma_i(C_i)$ is the graph of a $1$-Lipschitz
map on $\pi_{\leq d_i}(\sigma_i(C_i))$, where $d_i = \dim(C_i)$, 
\Cref{rema:1-Lipschitz_map_induces_bijection} yields the upper bound
\begin{equation*}
\#(\vox_{B^n}(\sigma(C_i)) \cap \vox_{B^n}(A^n))
\leq \#\bar \pi_{\leq d_i}(\vox_{B^n}(A^n))
= (\#\vox_B(A))^{d_i}.
\end{equation*}
Since each $d_i$ is bounded by $\dim(X)$, we obtain
\begin{equation*}
\#(\vox_{B^n}(X) \cap \vox_{B^n}(A^n))
= \sum_{i=1}^N \#(\vox_{B^n}(C_i) \cap \vox_{B^n}(A^n))
\leq N \cdot (\#\vox_B(A))^{\dim(X)}. \qedhere
\end{equation*}
\end{proof}

For the rest of this section, 
let $P,P',Q \subset K^n$ be cuboids around $0$
and let $B,B' \subset K^n$ be balls around $0$.
For a subset $X \subseteq K$ and tuple $t \in \{0,1\}^n$,
we write $X^t$ for $\prod_{i=1}^n t_iX$,
where $0 \cdot X = \{0\}$ and $1 \cdot X = X$.

\begin{lemm}\label{lemm:existence_of_inbetween_ball}
If $Q \subset P$, then there exists a ball $I \subset K$
such that $I^{t(P/Q)} \nsubseteq Q$ and $I^{t(P/Q)} \subseteq P$. 
\end{lemm}

\begin{proof}
Let $t = t(P/Q)$ and $D = \{ i \in \{1,\ldots,n\} \mid t_i = 1 \}$.
Since $Q$ is a proper subset of $P$, $D$ is non-empty.
Choose $j \in D$ such that the projection $\pi_j(P)$ is inclusion-minimal
among $\pi_i(P)$ for all $i \in D$. Set $I = \pi_j(P)$.
Then clearly, $I^t$ is contained in $P$.
Choose some $a \in \pi_j(P) \sm \pi_j(Q)$
and consider the point $x \in K^n$ with coordinates
\begin{equation*}
x_i = \begin{cases}
0, &\text{if $t_i = 0$}, \\
a_i, &\text{if $t_i = 1$}.
\end{cases}
\end{equation*}
Then $x \in I^t$ and $x \notin Q$.
\end{proof}

For a given sequence $(b_m)_{m \in \bb N}$, we write $\{ b_{<m} \}$ 
for $\{ b_{i} \mid 0 \leq i < m \}$.

\begin{lemm}\label{lemm:existence_of_good_sequence}
Assume that $K$ has infinite residue field.
If $\dim(P/B) > 0$, then there is a sequence $(b_m)_{m \in \bb N}$
in $K$, starting at $0$, such that $\{ b_{<m} \}^{t(P/B)} \subseteq P$ 
and $\#\vox_B(\{ b_{<m} \}^{t(P/B)}) = m^{\dim(P/B)}$ for all $m > 0$.
\end{lemm}

\begin{proof}
If $\dim(P/B) > 0$,
then $P$ is not a subset of $B$.
In turn, the intersection $B \cap P$ is a proper subset of $P$.
Hence, \Cref{lemm:existence_of_inbetween_ball} yields a ball $I \subset K$
such that $I^{t(P/B \cap P)} \nsubseteq B \cap P$ and $I^{t(P/B \cap P)} \subseteq P$.
In particular, there is some index $i$ such that $t(P/B \cap P)_i = 1$
and $\pi_i(B \cap P) \subsetneq I$.
Since $\pi_i(I) \subseteq \pi_i(P)$,
this leaves $\pi_i(B) \subsetneq I$.
Since $B$ is a ball, we have $B \subsetneq I^n$.

Since $K$ has infinite residue field, 
the quotient $I^n/B$ is infinite.
Hence, there is a sequence $(b_m)_{m \in \bb N}$ in $I$, starting at $0$,
such that each $b_m$ lies in a unique coset of $\pi_1(B)$.
Since $t(P/B \cap P) = t(P/B)$, 
the sequence $(b_m)_{m \in \bb N}$ satisfies the necessary conditions.
\end{proof}

\begin{prop}\label{prop:infinite_residue_field_isometric_cuboids_same_dim}
Assume that $\Th(K)$ is $1$-h-minimal and that $K$ has infinite residue field.
If $P$ and $P'$ are definably isometric, then $\dim(P/B) = \dim(P'/B)$.
\end{prop}

\begin{proof}
Without loss of generality, assume $\dim(P/B) \leq \dim(P'/B)$.
Let $f \colon P \to P'$ be a definable isometry.
By \Cref{rema:dim_existence_of_definable_cover}, 
the definable set $X = K^{t(P/B)} \cap P$ is of dimension $\dim(P/B)$
and $\vox_B(X)$ covers $\vox_B(P)$.
Since $f$ is an isometry, the image $f(X)$ is also of dimension $\dim(P/B)$,
and by the proof of \Cref{lemm:isometry_descends_to_bijection}, $\vox_B(f(X))$ covers $\vox_B(P')$.

If $\dim(P'/B) = 0$, then $\dim(P/B) = 0$ and we are done.
If $\dim(P'/B) > 0$, then there exists a sequence $(b_m)_{m \in \bb N}$ in $K$
as described in \Cref{lemm:existence_of_good_sequence}.
Then $\vox_B(\{b_{<m}\}^{t(P'/B)}) \subseteq \vox_B(f(X))$,
hence applying \Cref{prop:counting_balls} to $f(X)$ yields an integer $N \in \bb N$
such that
\begin{equation*}
\#\vox_B(\{b_{<m}\}^{t(P'/B)}) \leq N \cdot (\#\vox_{\pi_1(B)}(\{b_{<m}\}))^{\dim(f(X))}
\quad\text{for all $m > 0$.} 
\end{equation*}
By construction,
the left-hand side is equal to $m^{\dim(P'/B)}$
and the right-hand side is bounded by $N \cdot m^{\dim(f(X))}$.
Since $\dim(f(X))$ is equal to $\dim(P/B)$, we obtain
\begin{equation*}
m^{\dim(P'/B)} \leq N \cdot m^{\dim(P/B)}
\quad\text{for all $m > 0$.} 
\end{equation*}
Hence, $\dim(P'/B) \leq \dim(P/B)$. This finishes the proof.
\end{proof}

We can now prove our main result \Cref{theo:1.1}.

\begin{theo}\label{theo:infinite_residue_field_isometric_cuboids_equal_up_to}
Assume that $\Th(K)$ is $1$-h-minimal and let $P,P' \subset K^n$ be Cartesian products of balls in $K$.
If $P$ and $P'$ are definably isometric, then $P$ and $P'$ are equal 
up to translation and coordinate permutation.
\end{theo}

\begin{proof}
If $K$ has finite residue field, apply \Cref{prop:finite_residue_field_isometric_cuboids_equal_up_to}.
Assume that $K$ has infinite residue field and that $P$ and $P'$ are definably isometric.
As in the proof of \Cref{prop:finite_residue_field_isometric_cuboids_equal_up_to},
we may assume without loss of generality
that $P$ and $P'$ each contain $0$ and satisfy
$\pi_{i}(P) \subseteq \pi_{i+1}(P)$ and $\pi_{i}(P') \subseteq \pi_{i+1}(P')$ for $1 \leq i < n$.
Now suppose $P$ and $P'$ are not equal.
Then there is a largest index $i$ such that $\pi_i(P)$ and $\pi_i(P')$ differ. 
Without loss, assume $\pi_i(P) \subset \pi_i(P')$ and let $B = (\pi_i(P))^n$.
Then $\dim(P/B) = n-i$ and $\dim(P'/B) = n-(i-1)$. 
But this contradicts \Cref{prop:infinite_residue_field_isometric_cuboids_same_dim}.
\end{proof}

Since every ball in $K^n$ is the Cartesian product of its coordinate projections,
we further obtain the following generalization of \Cref{theo:1.1}.

\begin{coro}
Assume that $\Th(K)$ is $1$-h-minimal 
and let $B_1,\ldots,B_m \subset K^n$ and $B_1',\ldots,B_m' \subset K^n$
be finitely many balls.
If $\prod_{i=1}^m B_i$ and $\prod_{i=1}^m B_i'$ are definably isometric,
then $\prod_{i=1}^m B_i$ and $\prod_{i=1}^m B_i'$ are equal up to translation and coordinate permutation.
\end{coro}

\begin{proof}
The sets $\prod_{i=1}^m B_i$ and $\prod_{i=1}^m B_i'$ are cuboids in $K^{mn}$.
\end{proof}

Lastly, with a similar method as the one used in the proof of \Cref{prop:infinite_residue_field_isometric_cuboids_same_dim}, 
we obtain \Cref{prop:1.2}.

\begin{prop}\label{prop:dim_for_cuboids_relates_to_dim_for_definable_sets}
Assume that $\Th(K)$ is $1$-h-minimal and that $K$ has infinite residue field.
Then $\dim(P/B)$ is the minimal dimension of a definable set $X \subset K^n$
such that $\vox_B(X)$ covers $\vox_B(P)$.
\end{prop}

\begin{proof}
The existence of $X$ follows from \Cref{rema:dim_existence_of_definable_cover}.
Now let $X \subset K^n$ be an arbitrary definable set
such that $\vox_B(X)$ covers $\vox_B(P)$.
If $\dim(P/B) = 0$, then clearly $\dim(P/B) \leq \dim(X)$.
If $\dim(P/B) > 0$, then there is a sequence $(b_m)_{m \in \bb N}$ in $K$
as described in \Cref{lemm:existence_of_good_sequence}.
Applying \Cref{prop:counting_balls} to $X$ yields an integer $N \in \bb N$ such that
\begin{equation*}
\#\vox_B(\{b_{<m}\}^{t(P/B)}) 
\leq N \cdot (\vox_{\pi_1(B)}(\{b_{<m}\}))^{\dim(X)}
\quad\text{for all $m > 0$.} 
\end{equation*}
Since the left-hand side is equal to $m^{\dim(P/B)}$
and the right-hand side is bounded by $N \cdot m^{\dim(X)}$, we obtain
\begin{equation*}
m^{\dim(P/B)} 
\leq N \cdot m^{\dim(X)}
\quad\text{for all $m > 0$.} 
\end{equation*}
Hence, $\dim(P/B) \leq \dim(X)$. This finishes the proof.
\end{proof}